\documentclass{amsart}

\usepackage{fullpage}
\usepackage{amsmath}
\usepackage{amsfonts}
\usepackage{amssymb}

\input xy
\xyoption{matrix}\xyoption{arrow}\xyoption{curve}  \xyoption{frame}

\begin{document}

%%%%% Functors %%%%%

\newcommand{\Hom}{\mathrm{Hom}}
\newcommand{\RHom}{\mathrm{RHom}^*}
\newcommand{\HOM}{\mathrm{HOM}}
\newcommand{\stHom}{\underline{\mathrm{Hom}}}
\newcommand{\Ext}{\mathrm{Ext}}
\newcommand{\Tor}{\mathrm{Tor}}
\newcommand{\HH}{\mathrm{HH}}
\newcommand{\Endo}{\mathrm{End}}
\newcommand{\ENDO}{\mathrm{END}}
\newcommand{\stEnd}{\mathrm{\underline{End}}}
\newcommand{\Tr}{\mathrm{Tr}}

%%%%% Functions/Operators (no space after) %%%

\newcommand{\coker}{\mathrm{coker}}
\newcommand{\aut}{\mathrm{Aut}}
\newcommand{\op}{\mathrm{op}}
\newcommand{\add}{\mathrm{add}}
\newcommand{\ADD}{\mathrm{ADD}}
\newcommand{\ind}{\mathrm{ind}}
\newcommand{\rad}{\mathrm{rad}}
\newcommand{\soc}{\mathrm{soc}}
\newcommand{\ann}{\mathrm{ann}}
\newcommand{\im}{\mathrm{im}}
\newcommand{\chr}{\mathrm{char}}
\newcommand{\pdim}{\mathrm{p.dim}}

%%%%% Categories %%%%%%%%%

\newcommand{\rmod}{\mbox{mod-}}
\newcommand{\Rmod}{\mbox{Mod-}}
\newcommand{\lmod}{\mbox{-mod}}
\newcommand{\lMod}{\mbox{-Mod}}
\newcommand{\stmod}{\mbox{\underline{mod}-}}
\newcommand{\stlmod}{\mbox{-\underline{mod}}}

\newcommand{\gmod}[1]{\mbox{mod}_{#1}\mbox{-}}
\newcommand{\gMod}[1]{\mbox{Mod}_{#1}\mbox{-}}
\newcommand{\Bimod}[1]{\mathrm{Bimod}_{#1}\mbox{-}}

\newcommand{\proj}{\mbox{proj-}}
\newcommand{\lproj}{\mbox{-proj}}
\newcommand{\Proj}{\mbox{Proj-}}
\newcommand{\inj}{\mbox{inj-}}
\newcommand{\coh}{\mbox{coh-}}

%%%%% Style %%%%%%%%%%

\newcommand{\und}[1]{\underline{#1}}
\newcommand{\gen}[1]{\langle #1 \rangle}
\newcommand{\floor}[1]{\lfloor #1 \rfloor}
\newcommand{\ceil}[1]{\lceil #1 \rceil}
\newcommand{\bnc}[2]{\left(\scriptsize \begin{array}{c} #1 \\ #2 \end{array} \right)}
\newcommand{\bimo}[1]{{}_{#1}#1_{#1}}
\newcommand{\ses}[5]{\ensuremath{0 \rightarrow #1 \stackrel{#4}{\longrightarrow} 
#2 \stackrel{#5}{\longrightarrow} #3 \rightarrow 0}}
\newcommand{\A}{\mathcal{A}}
\newcommand{\B}{\mathcal{B}}
\newcommand{\uB}{\underline{\mathcal{B}}}
\newcommand{\C}{\mathcal{C}}
\newcommand{\D}{\mathcal{D}}
\newcommand{\tC}{\tilde{\mathcal{C}}}
\newcommand{\tK}{\tilde{K}(A)}
\newcommand{\ul}[1]{\underline{#1}}

%%%%% Theorems %%%%%%%%%%

\newtheorem{therm}{Theorem}[section]
\newtheorem{defin}[therm]{Definition}
\newtheorem{propos}[therm]{Proposition}
\newtheorem{lemma}[therm]{Lemma}
\newtheorem{coro}[therm]{Corollary}

\title{A construction of derived equivalent pairs of symmetric algebras}
\author{Alex Dugas}
\address{Department of Mathematics, University of the Pacific, 3601 Pacific Ave, Stockton, CA 95211, USA}
\email{adugas@pacific.edu}


\begin{abstract} Recently, Hu and Xi have exhibited derived equivalent endomorphism rings arising from (relative) almost split sequences as well as AR-triangles in triangulated categories.  We present a broader class of triangles (in algebraic triangulated categories) for which the endomorphism rings of different terms are derived equivalent.  We then study applications involving $0$-Calabi-Yau triangulated categories.  In particular, applying our results in the category of perfect complexes over a symmetric algebra gives a nice way of producing pairs of derived equivalent symmetric algebras.  Included in the examples we work out are some of the algebras of dihedral type with two or three simple modules.  We also apply our results to stable categories of Cohen-Macaulay modules over odd-dimensional Gorenstein hypersurfaces having an isolated singularity.
\end{abstract}

\maketitle
 
\section{Introduction} 
\setcounter{equation}{0}

A recent paper of Hu and Xi establishes a remarkable connection between almost split sequences and derived equivalences.  The main result of that paper, in its simplest incarnation, states that if $$\ses{X}{M}{Y}{f}{g}$$ is an almost split sequence of $\Lambda$-modules, then $\Endo_{\Lambda}(X \oplus M)$ is derived equivalent to $\Endo_{\Lambda}(Y \oplus M)$ \cite{HuXi1}.  In fact, this derived equivalence is provided by a tilting module  $T = \Hom_{\Lambda}(X \oplus M, Y \oplus M)$ of projective dimension $1$.  A similar pattern appears in the work of Iyama and Reiten on $3$-Calabi-Yau algebras \cite{IR}.  For an isolated Gorenstein singularity $R$ of dimension $3$, two noncommutative crepant resolutions $\Endo_R(M)$ and $\Endo_R(N)$ of $R$ are shown to be derived equivalent via the tilting module $\Hom_R(M,N)$ of projective dimension $1$.  In the case where $M$ and $N$ differ in only one indecomposable summand, these summands are related via a $2$-almost-split sequence (inside an appropriate subcategory) of reflexive $R$-modules.  %More precisely, if $M = M' \oplus X$, and $\ses{Y}{M_1 \longrightarrow M_0}{X}{}{f}$ is the $2$-almost split sequence ending in $X$ in $\add(M)$, then $N \cong M' \oplus \ker(f)$.

A common theme in these examples is the process of ``exchanging'' a summand of a module $T$ to obtain a new module whose endomorphism ring is derived equivalent to that of $T$.  However, the resulting derived equivalences, as well as those arising in the most general versions of Hu and Xi's theorem \cite{HuXi1} are all given by tilting modules.  In particular, they will never yield nontrivial derived equivalences between self-injective algebras.  Our goal in the present article is to show how these ideas can be extended to yield derived equivalences produced by tilting complexes.  Namely, given an algebraic triangulated category $\C$ containing a triangle of the form $$X \stackrel{f}{\longrightarrow} M \stackrel{g}{\longrightarrow} Y \longrightarrow X[1]$$ where every map $u : X \rightarrow M[i]$ factors through $f$ and every map $v: M[i] \rightarrow Y$ factors through $g$ (for $i \in \mathbb{Z}$), we produce two pairs of derived equivalent algebras
\begin{itemize}
\item $\Endo_{\C}(M\oplus X)$ and $\Endo_{\C}(M \oplus Y)$;
\item $\Endo_{\tC}(M\oplus X)$ and $\Endo_{\tC}(M \oplus Y)$.
\end{itemize}
Here $\tC$ denotes the orbit category $\C/[1]$ -- in particular, $\tC(X,Y) = \oplus_{i \in \mathbb{Z}} \C(X,Y[i])$ for all $X, Y \in \C$.  While the hypotheses may appear a bit forbidding at first glance, they are satisfied quite easily in some natural situations, including the category $K^b(\proj A)$ of perfect complexes over a symmetric algebra $A$ as well as the stable categories of Cohen-Macaulay modules over certain Gorenstein rings.  As an application, we illustrate a couple methods for obtaining interesting families of derived equivalent symmetric algebras.  We point out that our main results are quite similar in spirit to ongoing work of Hu, Koenig and Xi (see Remark (4) in \S 4 below), to whom we are very grateful for sharing a preprint of their article \cite{HuKX}.

This paper is organized as follows.  We start by reviewing the definitions of tilting modules and tilting complexes, and the construction of Okuyama-Rickard complexes in Section 2.  We also state Hu and Xi's theorem, which serves as our motivation and a model for our main results.  Section 3 is devoted to proving several lemmas for algebraic triangulated categories, which are necessary to deal with the non-uniqueness of the morphism of triangles that is guaranteed to exist by axiom (TR3) for triangulated categories.  We then state and prove our main result (Theorem 4.1) in Section 4.  In Section 5, we examine in more detail several settings where the hypotheses of Theorem 4.1 can be simplified, and we state the corresponding results (Theorem 5.2 and Corollary 5.4). The best cases are perhaps where $\C$ is $0$-Calabi-Yau, and we obtain pairs of derived equivalent symmetric algebras as a result.  We compute several examples in Section 6, chosen to illustrate some interesting derived equivalences between symmetric algebras that arise in this way.  Finally, we conclude with a brief discussion of derived equivalent endomorphism rings of perfect complexes, as produced in Theorem 5.2, from the point of view of dg-algebras.  

\vspace{2mm}
Throughout this article, we assume that $k$ is a field, and we consider primarily finite-dimensional algebras over $k$.  We typically work with right modules, unless noted otherwise.  In this case morphisms are written on the left and composed from right to left.  We also follow this convention for composition of morphisms in abstract categories, as well as for paths in quivers.  In particular, (covariant) representations of a quiver are identified with left modules over the path algebra.  We shall write $D$ for the duality $\Hom_k(-,k)$ on the category of finite-dimensional modules over a $k$-algebra.  For a category $\C$ and objects $X,Y \in \C$, we will write $\C(X,Y)$ for the set of morphisms from $X$ to $Y$.  All categories and functors we consider will be assumed to be $k$-linear over a field $k$, and subcategories will be assumed to be full and strict (i.e., closed uner isomorphisms).  Concerning complexes, we work with cochain complexes, i.e., the differential has degree $1$.  When we write out complexes, we will occasionally indicate the degree-$0$ term by underlining, and we frequently omit all terms which are zero.  We denote the morphism sets of a $k$-category $\A$ by $\A(X,Y)$ for objects $X,Y$ in $\A$.  We also write $K(\A)$ for the homotopy category of complexes in $\A$.  When $\A$ is abelian, we write $D(\A)$ for the derived category of $\A$.  If $A$ is a $k$-algebra we often substitute $K(A)$ for $K(\proj A)$ and $D(A)$ for $D(\rmod A)$.

\section{Background on tilting}
\setcounter{equation}{0}

  We begin by briefly recalling the definitions of tilting modules and complexes and their connection with derived equivalence.  We then review a theorem of Hu and Xi that produces tilting modules and hence examples of derived-equivalent pairs of algebras.  For simplicity, we assume here that $A$ is a $k$-algebra.  Recall that an $A$-module $T_A$ is a {\it tilting module} if (1) $\pdim\ T_A < \infty$; (2) $\Ext^i_{A}(T,T) = 0$ for all $i \geq 1$; and (3) there exists an exact sequence $0 \rightarrow A \longrightarrow T_0 \longrightarrow \cdots \longrightarrow T_r \rightarrow 0$ with all $T_i \in \add(T)$.  A well-known result of Happel's states that the (bounded) derived categories $D^b(A)$ and $D^b(\Endo_A(T))$ are triangle equivalent for any tilting module $T_A$ \cite{TCRTA}.  More generally, Rickard has shown that two algebras $A$ and $B$ are derived equivalent if and only if $B \cong \Endo_{K(A)}(T^{\bullet})$ where $T^{\bullet} \in K^b(\proj A)$ is a {\it tilting complex}, meaning that (1) $\Hom_{K^b(A)}(T,T[n]) = 0$ for all $n \neq 0$; and (2) $T$ generates $K^b(\proj A)$ as a triangulated category \cite{MTDC}.  One can easily check that the projective resolution of a tilting module satisfies this definition.  However, over self-injective algebras one must look elsewhere for tilting complexes, since such algebras admit no nonprojective modules of finite projective dimension.

  Furthermore, recall that for a subcategory $\C$ of an additive category $\A$, a map $f : X \rightarrow C$ (resp. $g : C \rightarrow X$) is a {\it left} (resp. {\it right}) {\it $\C$-approximation} of $X$ if $C \in \C$ and the induced map
   $$\A(C,-) \stackrel{(f,-)}{\longrightarrow} \A(X,-)\ \ \mbox{(resp.\ }\A(-,C) \stackrel{(-,g)}{\longrightarrow} \A(-,X)\ )$$ is a surjective morphism of functors on $\C$.

The following proposition is probably well-known and provides a simple construction of tilting complexes of length $1$ over self-injective algebras.  Complexes of this form have been used extensively to verify various cases of Brou\'e's abelian defect group conjecture \cite{RFGT}, as well as in realizing derived equivalences between various symmetric algebras \cite{Holm1, Holm2}.  As these particular complexes will play a central role in our main results, we supply a proof.  Note, in particular, that if $A$ is (weakly) symmetric, then the hypothesis $\nu P \cong P$ is automatically satisfied by any projective $P$.

\begin{propos}[Cf. \cite{RFGT}]  Suppose that $A$ is a basic self-injective algebra with Nakayama functor $\nu$.  If $A_A = P \oplus Q$ with $\nu P \cong P$, and $P \stackrel{f}{\longrightarrow} Q'$ is a left $\add(Q)$-approximation, then $$T = [P \stackrel{f}{\longrightarrow} \und{Q'}] \oplus \und{Q}$$ is a tilting complex.
\end{propos}

\noindent
{\it Proof.}  We set $T_1 = [P \stackrel{f}{\longrightarrow} \und{Q'}]$ and $T_2 = \und{Q}$.  It is clear that $T$ generates $K^b(\proj A)$ as a triangulated category, so we must simply check that $\Hom_{K(A)}(T_1,T[1]) = 0$ and $\Hom_{K(A)}(T[1],T_1) = 0$.  The former follows since any map $P \rightarrow Q$ factors through $f$ by hypothesis.  To see the latter, let $g : Q \rightarrow P$ be a nonzero map such that $fg = 0$.  Then there is a map $h : P \rightarrow \nu Q \cong Q$ such that $hg \neq 0$.  But such an $h$ factors through $f$ by hypothesis, contradicting $fg = 0$.  $\Box$\\

We now state Hu and Xi's theorem.

\begin{therm}[Hu, Xi \cite{HuXi1}] Suppose $\ses{X}{M'}{Y}{f}{g}$ is an almost $\add(M)$-split sequence in $\rmod A$, i.e.,
\begin{itemize}
\item[(i)] $f$ is a left $\add(M)$-approximation, meaning $\Hom_A(f,M)$ is onto; and
\item[(ii)] $g$ is a right $\add(M)$-approximation, meaning $\Hom_A(M,g)$ is onto.
\end{itemize}
Then $\Endo_A(M \oplus X)$ and $\Endo_A(M \oplus Y)$ are derived equivalent.  
\end{therm}

In fact, Hu and Xi show that this derived equivalence is afforded by a tilting module $T$ with $\pdim\ T \leq 1$, which is defined by the following projective presentation in $\rmod \Endo_A(M \oplus X)$:
$$\ses{\Hom_A(M \oplus X, X)}{\Hom_A(M \oplus X, M' \oplus M)}{T}{\bnc{(M \oplus X, f)}{0}}{}.$$
Furthermore, observe that the condition that $f$ is a left $\add(M)$-approximation corresponds to the map $(M \oplus X, f)$ being a left $\add((M \oplus X, M))$-approximation of $(M \oplus X, X)$ in the category of (projective) right $\Endo_A(M \oplus X)$-modules.  Thus the tilting module $T$ is obtained from Riedtmann and Schofield's construction: replacing the summand $(M \oplus X, X)$ of $\Endo_A(M \oplus X)$ with a non-projective summand.  In case $X$ is indecomposable, $T$ will be the minimal nonprojective completion of the almost complete tilting module $(M \oplus X,M)$ over $\Endo_A(M \oplus X)$.  
  
The most natural setting in which this theorem applies is when the short exact sequence above is an almost split sequences in $\rmod A$.  Hu and Xi also obtain an analagous conclusion when $X \stackrel{f}{\longrightarrow} M \stackrel{g}{\longrightarrow} Y \stackrel{w}{\longrightarrow} X[1]$ is an AR-triangle in a Hom-finite, triangulated Krull-Schmidt category such that $X[1] \not \in \add(M\oplus Y)$.  Additionally, Hu and Xi find examples of tilting modules of projective dimension $n \geq 2$ by considering $n$-almost split sequences (as introduced by Iyama \cite{Iyama1}), but we will not pursue generalizations along these lines here.

\section{Technical necessities for triangulated categories}
\setcounter{equation}{0}

We start by reviewing some properties of (algbebraic) triangulated categories and prove a lemma that will be essential in the proof of our main result.  We let $(\mathcal{T},\Sigma)$ be a triangulated category with suspension functor $\Sigma$.  Recall that $\mathcal{T}$ is said to be {\it algebraic} if it is triangle equivalent to the stable category of an exact Frobenius category $(\mathcal{B},S)$, as described in Happel's book \cite{TCRTA} for instance.  Following Happel's notation, $\mathcal{B}$ denotes an extension closed full subcategory of an abelian category, and $S$ denotes the set of short exact sequences (in this abelian category) all of whose terms belong to $\mathcal{B}$.  That $(\mathcal{B},S)$ is Frobenius means that $\mathcal{B}$ has enough $S$-projectives and enough $S$-injectives and that these objects coincide.  The (injective) stable category of $(\mathcal{B},S)$, obtained by factoring out the morphisms in $\mathcal{B}$ that factor through an injective, will be denoted $\ul{\mathcal{B}}$, and we write $\ul{f} \in \uB(X,Y)$ for the image of a morphism $f \in \B(X,Y)$.  In particular, whenever we specify a morphism $\ul{f}$ of $\uB$, we are implicitly choosing a representative morphism $f$ of $\B$.  The suspension functor $\Sigma$ of $\uB$ is given by the cosyzygy functor $\Omega^{-1}$.  Note that its definition requires us to fix sequences $\ses{X}{I(X)}{\Omega^{-1}X}{\mu_X}{\pi_X}$  for each object $X$ of $\mathcal{B}$, although these choices do not affect the isomorphism type of $\Sigma$.   Furthermore, Happel shows that these choices can be made to guarantee that $\Sigma = \Omega^{-1}$ is an automorphism of $\ul{\mathcal{B}}$, under the assumption that there exist bijections between the isomorphism classes in $\ul{\mathcal{B}}$ of $X$ and $\Omega^{-1}X$ for any $X \in \mathcal{B}$.  The latter assumption will clearly be satisfied inside a basic category (i.e., where each isomorphism class consists of one object), and so by passing to a skeleton of $\ul{\mathcal{B}}$, if necessary, we may assume that $\Sigma = \Omega^{-1}$ is an automorphism of $\ul{\mathcal{B}}$.  We do not need this right away, but it is important for our applications in the next section.

  Recall that the standard triangles in $\ul{\mathcal{B}}$ are defined by completing a morphism $u : X \rightarrow Y$ in $\mathcal{B}$ via a pushout diagram in $\mathcal{B}$ as shown below.
\begin{eqnarray}\vcenter{
\xymatrix{0 \ar[r] & X \ar[r]^{\mu_X} \ar[d]^u & I(X) \ar[r]^{\pi_X} \ar[d]^t & \Sigma X \ar@{=}[d] \ar[r] & 0 \\
0 \ar[r] & Y \ar[r]_v & C_u \ar[r]_w & \Sigma X \ar[r] & 0}}
\end{eqnarray}
Of course, the distinguished triangles are then defined to be those that are isomorphic to a standard triangle $X \stackrel{\ul{u}}{\longrightarrow} Y \stackrel{\ul{v}}{\longrightarrow} C_u \stackrel{\ul{w}}{\longrightarrow} \Sigma X \rightarrow$.

Our main concern in this section is related to the axiom (TR3) of triangulated categories, which states that any commutative square $(*)$ in $\uB$ can be completed to a morphism of triangles.
\begin{eqnarray}\vcenter{
\xymatrix{ X \ar[r]^{\ul{u}} \ar[d]^{\ul{f}} \save+<5ex,-4ex> \drop{(*)} \restore & Y \ar[r]^{\ul{v}} \ar[d]^{\ul{g}} & Z \ar[r]^{\ul{w}} \ar@{-->}[d]^{\ul{h}} & \Sigma X \ar[d]^{\Sigma \ul{f}} \\
X' \ar[r]^{\ul{u'}} & Y' \ar[r]^{\ul{v'}} & Z' \ar[r]^{\ul{w'}} & \Sigma X' } }
\end{eqnarray}
Such a completion is usually not unique, and resolving the ambiguity that thus arises is the focus of an article by Neeman \cite{NATC}.  In particular, he points out that in the category $K(\A)$ of chain complexes over an additive category $\A$ with homotopy classes of chain maps, the standard mapping cone construction leads to a natural set of choices for the third map $h$, which is closed under addition and composition.  The construction of these ``naturally good'' completions uses the fact that $K(\A)$ is a quotient of the category of chain complexes and chain maps; thus it is no surprise that the same idea can be adapted to any algebraic triangulated category as we now verify.
  
We start by reviewing the ``standard'' completion of a commutative square to a morphism of standard triangles in $\ul{\B}$ given by Happel in his verification of (TR3) \cite{TCRTA}. Following the notation of the above diagrams, a commutative square $(*)$ implies that $gu - u'f = \alpha \mu_X$ for some $\alpha : I(X) \rightarrow Y'$, which is not necessarily unique. (Note that we are working with morphisms in $\B$ for now).  We can also lift $f$ to a map $I_f : I(X) \rightarrow I(X')$ such that $I_f \mu = \mu' f$, although again this map need not be unique.  Note that $I_f$ induces a map $\Sigma f : \Sigma X \rightarrow \Sigma X'$, whose image in $\ul{\B}$ is independent of the choice of $I_f$, and hence denoted $\Sigma \ul{f}$.  Now the two maps $v'g : Y \rightarrow C_{u'}$ and $v' \alpha + t'I_f : I(X) \rightarrow C_{u'}$ induce a unique $h : C_u \rightarrow C_{u'}$ such that $v'g = hv$ and $v'\alpha +t'I_f = ht$, by the universal property of pushouts.  One can then apply this universal property once more to check that $w'h = (\Sigma f)w$, so that $(\ul{f}, \ul{g}, \ul{h})$ gives a morphism of triangles.  Alternatively, once we have chosen $\alpha$ and $I_f$, $h$ is the unique map making the following diagram commutative in $\B$.

\begin{eqnarray}
\vcenter{\xymatrixrowsep{2.0pc} \xymatrixcolsep{3.0pc} \xymatrix{0 \ar[r] & X \ar[r]^{\bnc{-\mu_X}{u}} \ar[d]^f & I(X) \oplus Y \ar[r]^{(t, v)} \ar[d]^{\bnc{I_f\ 0}{\alpha\ g}} & C_u \ar[r] \ar@{-->}[d]^h & 0 \\
0 \ar[r] & X' \ar[r]_{\bnc{-\mu_{X'}}{u'}} & I(X') \oplus Y' \ar[r]_{(t', v')} & C_{u'} \ar[r] & 0}}
\end{eqnarray}

We shall informally refer to the maps $h$ constructed in this way as {\it good maps} (relative to the pair $(f,g)$).  Using the above description of good maps via short exact sequences, it is easy to see that the sum of two good maps $h, h' : C_u \rightarrow C_{u'}$ is good (relative to the sum of the corresponding pairs), and the composition of good maps $h : C_u \rightarrow C_{u'}$ and $h': C_{u'} \rightarrow C_{u''}$ is good as well (relative to the composite of the corresponding pairs).  When we refer to a good map we will frequently omit the reference to the pair $(f,g)$ when there is no danger of confusion.  However, we point out that our notion of {\it good} is far from absolute, as many (often, all) maps $C_u \rightarrow C_{u'}$ may be good relative to different pairs $(f,g)$ (cf. Lemma 3.4).

Now recall that we had to make choices for $\alpha : I(X) \rightarrow Y'$ and $I_f : I(X) \rightarrow I(X')$.  In general, the good map $h$ (and even $\ul{h}$) constructed above depends on these choices.  If $\alpha'$ is another choice for $\alpha$, then $(\alpha - \alpha')\mu = 0$, whence $\alpha - \alpha' = \beta \pi$ for some map $\beta : \Sigma X \rightarrow Y'$.  Likewise, if $I'_f$ is another choice for $I_f$, then $(I_f - I'_f)\mu = 0$, whence $I_f - I'_f = \gamma \pi$ for some $\gamma : \Sigma X \rightarrow I(X')$.  Suppose that the choices of $\alpha'$ and $I'_f$ lead to the map $h' : C_u \rightarrow C_{u'}$ with $v'g = h'v$ and $v'\alpha' +t'I'_f = h't$.  Then
$$(h-h')t = v'(\alpha-\alpha')+t'(I_f - I'_f) = v'\beta \pi + t' \gamma \pi = (v'\beta w + t'\gamma w) t,$$ and $(h-h')v = 0 = (v'\beta w + t'\gamma w) v$.  Hence, by the universal property of pushouts, $h-h' = v'\beta w + t'\gamma w$ in $\B$, which yields $\ul{h} = \ul{h'} + \ul{v'}\ul{\beta}\ul{w}$ in $\ul{\B}$.  In other words, we see that the ambiguity in the construction of a good map $h : C_u \rightarrow C_{u'}$ in $\ul{\B}$ comes precisely from the maps of the form $\ul{v'}\ul{\beta}\ul{w}$ with $\beta : \Sigma X \rightarrow Y'$.  The following lemma is now immediate.

\begin{lemma} Suppose we have a commutative square $(*)$ in $\ul{\B}$ which we wish to extend to a morphism of triangles as in the diagram below.
$$\xymatrix{ X \ar[r]^{\ul{u}} \ar[d]^{\ul{f}} \save+<5ex,-4ex> \drop{(*)} \restore & Y \ar[r]^{\ul{v}} \ar[d]^{\ul{g}} & C_u \ar[r]^{\ul{w}} \ar@{-->}[d]^{\ul{h}} & \Sigma X \ar[d]^{\Sigma f} \\
X' \ar[r]^{\ul{u'}} & Y' \ar[r]^{\ul{v'}} & C_{u'} \ar[r]^{\ul{w'}} & \Sigma X' }$$
If $\ul{\B}(\ul{w},Y') = 0$, then there is a unique map $\ul{h}$ in $\uB$ making the diagram commute, for which $h$ is good relative to $(f,g)$.
\end{lemma}

%It is not hard to see that if two pairs of maps $(f,g)$ and $(f'g')$ induce maps $h$ and $h'$ by the above construction, then $(f+f',g+g')$ along with the choices $\alpha + \alpha'$ and $I_f + I_{f'}$ will yield $h+h'$.  In the following lemma, we check that this construction also behaves well under composition.

%\begin{lemma}  Suppose that the above construction yields two maps $h : C_u \rightarrow C_{u'}$ and $h' : C_{u'} \rightarrow C_{u''}$ completing morphisms of triangles as in the diagram below.  Then the same construction, starting with $(f'f, g'g)$ and making a suitable choice of $\alpha'' : I(X) \rightarrow Y''$, yields the map $h'h$.
%$$\xymatrix{ X \ar[r]^u \ar[d]^f & Y \ar[r]^v \ar[d]^g & C_u \ar[r]^w \ar@{-->}[d]^h & \Sigma X \ar[d]^{\Sigma f} \\ X' \ar[r]^{u'} \ar[d]^{f'} & Y' \ar[r]^{v'} \ar[d]^{g'} & C_{u'} \ar[r]^{w'} \ar[d]^{h'} & \Sigma X' \ar[d]^{\Sigma f'}  \\ X'' \ar[r]^{u''} & Y'' \ar[r]^{v''} & C_{u''} \ar[r]^{w''} & \Sigma X''}$$ \end{lemma}

%\noindent {\it Proof.}  Suppose $\alpha, \alpha', I_f$ and $I_{f'}$ have been chosen as above in the construction of $h$ and $h'$.  Since we have $g'gu - u''f'f = g'(gu - u'f) + (g'u' - u''f')f = g'\alpha \mu_X + \alpha'\mu_{X'}f = (g'\alpha + \alpha'I_f)\mu_X$, we can choose $\alpha'' := g'\alpha +\alpha'I_f$.  Furthermore, we may take $I_{f'f} = I_f'I_f$.  We now simply need to check that $h'hv = v''g'g$ and $h'ht = v''\alpha'' + t''I_{f'f}$.  The first is clear, while for the second we have $h'ht = h'(v'\alpha + t'I_f) = v''g'\alpha + v''\alpha'I_f + t''I_{f'} I_f = v''\alpha'' + t''I_{f'f}$. $\Box$\\

It follows that, under the hypothesis $\ul{\B}(\ul{w}, Y') = 0$, there is a well-defined mapping from the set of pairs $(f,g)$, for which $(\ul{f}, \ul{g})$ gives a chain map (over $\uB$) from $X \stackrel{\ul{u}}{\longrightarrow} Y$ to $X' \stackrel{\ul{u'}}{\longrightarrow} Y'$, into $\ul{\B}(C_u,C_{u'})$.  The next two lemmas address the kernel of this map.

\begin{lemma} Suppose that the commutative square $(*)$ in $\und{\B}$ is null-homotopic in the sense that there is a map $\varphi : Y \rightarrow X'$ with $\ul{f} = \ul{\varphi u}$ and $\ul{g} = \ul{u'\varphi}$.  Then $\alpha$ can be chosen so that the above construction yields a good map $h$, relative to $(f,g)$, with $\ul{h} = 0$.
\end{lemma}

\noindent
{\it Proof.}  By assumption we can write $g = u'\varphi + \beta \mu_Y$ and $f = \varphi u + \gamma \mu_X$ for certain maps $\beta : I(Y) \rightarrow Y'$ and $\gamma : I(X) \rightarrow X'$.  In addition, we can lift $u$ to a map $I_u : I(X) \rightarrow I(Y)$ such that $\mu_Y u = I_u \mu_X$.  Then $$gu - u'f = \beta \mu_Yu - u'\gamma \mu_X = (\beta I_u - u'\gamma)\mu_X,$$ showing that we may take $\alpha := \beta I_u - u'\gamma$.  Now notice that since $\mu_Y u = I_u \mu_X$, the pushout property gives a map $r : C_u \rightarrow I(Y)$ such that $rv = \mu_Y$ and $rt = I_u$.  Similarly, since $(I_f - \mu_{X'} \gamma)\mu_X = \mu_{X'} \varphi u$, we obtain a map $s : C_u \rightarrow I(X')$ such that $st = I_f - \mu_{X'}\gamma$ and $sv = \mu_{X'}\varphi$.  We now claim that $h = v'\beta r + t's : C_u \rightarrow C_{u'}$, which is clearly zero in $\ul{\B}$.  To see the claim, observe $v'g = v'u'\varphi + v'\beta \mu_Y = t' \mu_{X'} \varphi + v'\beta r v = (t's + v'\beta r)v$ and $v'\alpha + t'I_f = v'\beta I_u - t' \mu_{X'} \gamma + t'I_f = (v'\beta r + t's)t$.  $\Box$\\

In particular, under the hypothesis of Lemma 3.1, if $\ul{f} = 0 = \ul{g}$, the unique map $\ul{h}$ must be $0$.  %From our previous remarks about good maps,
It now follows that, under the hypothesis of Lemma 3.1, the unique map $\ul{h}$ is determined by $\ul{f}$ and $\ul{g}$, independent of the choices of maps $f$ and $g$ in $\B$.  In the sequel, we may thus refer to this map $\ul{h}$ as the unique good map determined by $\ul{f}$ and $\ul{g}$.

\begin{lemma} If the commutative square $(*)$ can be completed by a good map $h : C_u \rightarrow C_{u'}$ with $\ul{h} = 0$, then the square $(*)$ is null-homotopic (as defined in the previous lemma).
\end{lemma}

\noindent
{\it Proof.}  Suppose some choice of $\alpha$ and $I_f$ produces a map $h$ with $\ul{h} =0$.  Then $h = h_1 \mu_C$ for some map $h_1 : I(C_u) \rightarrow C_{u'}$.  Since $I(C_u)$ is projective, we can further factor $h_1$ over the epimorphism $(t', v')$ and we write $h_1 = t's' + v'q'$ for maps $s': I(C_u) \rightarrow I(X')$ and $q': I(C_u) \rightarrow Y'$.  We set $s = s'\mu_C$ and $q = q' \mu_C$.  Now, since $(t', v')\bnc{s}{q}(t, v) = (t', v')\bnc{I_f\ 0}{\alpha\ g}$, we must have maps $\gamma: I(X) \rightarrow X'$ and $\varphi : Y \rightarrow X'$ with $$\bnc{I_f\ 0}{\alpha\ g} - \bnc{s}{q}(t, v) = \bnc{-\mu_{X'}}{u'}(\gamma, \varphi).$$  In particular, we find $qv = g - u'\varphi$, and thus $\ul{g} = \ul{u'\varphi}$ since $\ul{q}=0$.  We also have $$\bnc{-\mu_{X'}}{u'}f = \left[ \bnc{I_f\ 0}{\alpha\ g} - \bnc{s}{q}(t, v)\right] \bnc{-\mu_x}{u} = \bnc{-\mu_{X'}}{u'}(\gamma, \varphi) \bnc{-\mu_x}{u},$$ which implies that $f = \varphi u - \gamma \mu_X$ as $\bnc{-\mu_{X'}}{u'}$ is a monomorphism.   Thus $\ul{f} = \ul{\varphi u}$, and $\varphi$ is the required homotopy. $\Box$\\

\begin{lemma}  Suppose that $X \stackrel{\ul{u}}{\longrightarrow} Y \stackrel{\ul{v}}{\longrightarrow} C_u \stackrel{\ul{w}}{\longrightarrow}$ and $X' \stackrel{\ul{u'}}{\longrightarrow} Y' \stackrel{\ul{v'}}{\longrightarrow} C_{u'} \stackrel{\ul{w'}}{\longrightarrow}$ are standard triangles, and let $g : Y \rightarrow Y'$ and $h : C_u \rightarrow C_{u'}$ be maps in $\B$ for which $\ul{v'g} = \ul{hv}$.  Then there exist maps $f: X \rightarrow X'$ and $g' : Y \rightarrow Y'$ and $h' : C_u \rightarrow C_{u'}$ such that $h'$ is good relative to the pair $(f,g')$ and $\ul{h'} = \ul{h}$.
\end{lemma}

\noindent
{\it Proof.}  Since $hv - v'g$ factors through a projective, it factors through the epimorphism $(t',v') : I(X') \oplus Y' \rightarrow C_{u'}$ in (3.3).  Thus we may write $hv -v'g = t'r + v's$ for $\bnc{r}{s} : Y \rightarrow I(X') \oplus Y'$.  Now set $g' = g+s$ to get $hv = v'g' + t'r$.  From (3.1) we see that $v$ is a monomorphism, and hence we can write $r = r'v$ for $r' : C_u \rightarrow I(X')$ since $I(X')$ is injective.  Now set $h' = h - t'r'$ and notice $\ul{h'} = \ul{h}$, while we now have $h'v = v'g'$.  Furthermore, as a map from a projective object $h't$ factors through the epimorphism $(t',v)$ via a map $\bnc{\phi}{\alpha} : I(X) \rightarrow I(X') \oplus Y'$.  This makes the right square in (3.3) commutative (with $g'$ and $h'$ in place of $g$ and $h$ respectively), and hence we obtain an induced map $f: X \rightarrow X'$, making the left square of (3.3) commute.  Clearly, $h'$ is good relative to the pair $(f,g')$. $\Box$\\

In the following section, we will use these lemmas to obtain a ring isomorphism between the homotopy-classes of endomorphisms of certain complexes of the form $X \stackrel{\ul{u}}{\rightarrow} Y$ in $\ul{\B}$ and the endomorphism ring of $C_u$ in $\ul{\B}$.

\section{Two variations on Hu and Xi's Theorem}
\setcounter{equation}{0}

In this section we formulate and prove two variations of Hu and Xi's theorem that yield derived equivalences furnished by tilting complexes of length $1$.  Let $\mathcal{C}$ be an algebraic Krull-Schmidt triangulated category with suspension now denoted $[1]$.  As we no longer need to represent $\C$ as a stable category, morphisms in $\C$ will not be underlined.  We further assume that the suspension is an automorphism of $\mathcal{C}$ (see the comments at the beginning of the previous section).  We can now define $\tC$ to be the {\it orbit category} $\mathcal{C}/[1]$ defined by ``factoring out'' the action of the suspension; i.e., $\tC$ has the same objects as $\mathcal{C}$, but the morphism sets are $$\tC(X,Y) = \oplus_{i \in \mathbb{Z}} \C(X,Y[i]).$$
We let $\pi: \mathcal{C} \rightarrow \tC$ be the natural covering functor, which is the identity on objects and coincides with the obvious inclusion $\mathcal{C}(X,Y) \hookrightarrow \tC(X,Y)$ on morphisms.  When there is no risk of confusion, we will typically omit $\pi$ from our notation when considering the images of objects or morphsims of $\mathcal{C}$ inside $\tC$.  Note that $\tC$ is a $\mathbb{Z}$-graded category and any endomorphism ring $\tC(X,X)$ becomes a $\mathbb{Z}$-graded $k$-algebra.  We also remind the reader that if $X$ and $Y$ are indecomposable objects of $\C$, then $X \cong Y$ in $\tC$ if and only if $X \cong Y[i]$ for some $i \in \mathbb{Z}$.  

To state our theorem, we will need a stronger version of left and right $\D$-approximations for a subcategory $\D$ of $\C$.  Letting $\gen{\D}$ denote the full, additive subcategory generated by $\cup_{i \in \mathbb{Z}}\D[i]$, we will rely on left and right $\gen{\D}$-approximations instead.  In particular, a map $f : X \rightarrow D$ is a left $\gen{\D}$-approximation if $D \in \gen{\D}$ and every map $g : X \rightarrow D'[i]$ with $D' \in \D$ factors through $f$.  In terms of the orbit category $\tC$, we write $\pi(\D)$ for the strict, full subcategory of $\tC$ generated by the objects of $\D$, and then it can be seen that $f$ is a left $\gen{\D}$-approximation if and only if $\pi(f)$ is a left $\pi(\D)$-approximation of $X$ in $\tC$.  The situation is similar for right $\gen{\D}$-approximations and right $\pi(\D)$-approximations.  We point out that a $\gen{\D}$-approximation is somewhat different than a  {\it cohomological $\D$-approximation} as used in \cite{HuKX}.  For one, the approximating object $D$, in our case, need only belong to $\gen{\D}$, as opposed to $\D$. %In both situations, the advantage of considering these `stronger' approximations is given by the exactness of the functors $\C(-,M[i])$ and $\C(M[i],-) on an approximating triangle. 

\begin{therm}  Suppose $\mathcal{C}$ contains a triangle $$X \stackrel{f}{\longrightarrow} M' \stackrel{g}{\longrightarrow} Y \stackrel{h}{\longrightarrow} X[1]$$
where $M' \in \gen{M}$ for some $M \in \C$, and
\begin{itemize} 
\item[(a)] $f$ is a left $\gen{M}$-approximation; and 
\item[(b)] $g$ is a right $\gen{M}$-approximation.
\end{itemize}
  Then \begin{enumerate}
\item $\Lambda = \Endo_{\tC}(M \oplus X)$ and $\Gamma = \Endo_{\tC}(M \oplus Y)$ are derived equivalent.  More specifically,
\begin{enumerate}
\item[(i)] $T = \ul{\tC(M \oplus X, M)} \oplus [\tC(M \oplus X, X) \stackrel{\pi(f)_*}{\longrightarrow} \ul{\tC(M \oplus X, M')}] =: T_1 \oplus T_2 $ is a tilting complex over $\Lambda$, concentrated in degrees $-1$ and $0$; and
\item[(ii)] $\Endo_{K(\Lambda)}(T) \cong \Gamma$.
\end{enumerate}
\vspace{3mm}
\item For any $M'' \in \gen{M}$ with $M' \in \add(M'')$, $\Lambda' = \Endo_{\C}(M'' \oplus X)$ and $\Gamma' = \Endo_{\C}(M'' \oplus Y)$ are derived equivalent.  More specifically,
\begin{enumerate}
\item[(i)] $T' = \ul{\C(M'' \oplus X, M'')} \oplus [\C(M'' \oplus X, X) \stackrel{f_*}{\longrightarrow} \ul{\C(M'' \oplus X, M')}] =: T'_1 \oplus T'_2 $ is a tilting complex over $\Lambda'$, concentrated in degrees $-1$ and $0$; and
\item[(ii)] $\Endo_{K(\Lambda')}(T') \cong \Gamma'$.
\end{enumerate} \end{enumerate}
Furthermore, the hypothesis (b) (respectively, (a)) can be deduced from (a) (resp., (b)) in case $\Lambda$ (resp., $\Gamma$) is weakly symmetric.
\end{therm}

\noindent
{\it Proof.}  Notice that (a) is equivalent to $\C(f,M[i]): \C(M',M[i]) \rightarrow \C(X,M[i])$ being an epimorphism for all $i \in \mathbb{Z}$.  Equivalently, the long exact sequence resulting from applying $\C(-,M)$ to the given triangle splits up into the short exact sequences isomorphic to 
\begin{eqnarray} \ses{\C(Y,M[i])}{\C(M',M[i])}{\C(X,M[i])}{\C(g,M[i])}{\C(f,M[i])}
\end{eqnarray}
 for all $i \in \mathbb{Z}$.  Similarly, (b) is equivalent to the existence of short exact sequences 
\begin{eqnarray} \ses{\C(M[i],X)}{\C(M[i],M')}{\C(M[i],Y)}{\C(M[i],f)}{\C(M[i],g)}
\end{eqnarray}
 for all $i \in \mathbb{Z}$.  Taking the direct sums of these short exact sequences over all $i \in \mathbb{Z}$ now shows that (a) is equivalent to the exactness of
\begin{eqnarray} \ses{\tC(Y,M)}{\tC(M',M)}{\tC(X,M)}{\tC(\pi(g),M)}{\tC(\pi(f),M)} \end{eqnarray} and (b) is equivalent to the exactness of \begin{eqnarray} \ses{\tC(M,X)}{\tC(M,M')}{\tC(M,Y)}{\tC(M,\pi(f))}{\tC(M,\pi(g))}. \end{eqnarray}
If we know that $\Lambda$ is weakly symmetric, then we have an isomorphism of functors $D\tC(-,M) \cong \tC(M,-)$ on $\add(M \oplus X)$.  Thus, applying the duality $D$ to (4.3), we see that $\tC(M,f)$ is monic, whence $\C(M[i],g)$ is an epimorphism for all $i$ and (b) follows.  Similarly, if we assume (b) and that $\Gamma$ is weakly symmetric, then (a) would follow automatically.

We now prove (1).  For (i), first note that $T$ generates $K^b(\proj \Lambda)$ since $\tC(M \oplus X, X)[1]$ can be recovered as the mapping cone of a map from $\tC(M \oplus X, M') \in \add(T_1)$ to $T_2$.  Next, notice that any map $\alpha : T_2 \rightarrow T_1[1]$ will be induced by a map in $\tC(X,M)$, i.e., by a sum of maps in $\C(X,M[i])$ for various $i$.  But such a map factors through $f$ by (a), and hence $\alpha$ will be null-homotopic.  Similarly, we see that $\Hom_{K(\Lambda)}(T_2,T_2[1])=0$.  Likewise, any map $\beta : T[1] \rightarrow T_2$ will be induced by a map $\beta \in \tC(M, X)$ such that $ \pi(f) \beta= 0$.  However, (4.4) shows that such a $\beta$ must be zero.

We next focus on describing the endomorphism ring of $T$ in order to verify (ii).  Corresponding to the decomposition $T = T_1 \oplus T_2$, we can express this endomorphism ring in matrix form:
\begin{eqnarray*}
\Endo_{K(\Lambda)}(T) \cong \left( \begin{array}{rr} \Endo_{\Lambda}(T_1) & \Hom_{K(\Lambda)}(T_2,T_1) \\ \Hom_{K(\Lambda)}(T_1,T_2) & \Endo_{K(\Lambda)}(T_2) \end{array} \right). 
\end{eqnarray*}
We begin by describing each entry in the above matrix.  First, notice that $\Endo_{\Lambda}(T_1) \cong \Endo_{\tC}(M)$ as rings.  Next, observe that any map $\alpha$ in $K(\Lambda)$ from $T_1$ to $T_2$ consists of a map from $\tC(M\oplus X, M)$ to $\tC(M \oplus X, M')$, and such a map is induced by a unique $u \in \tC(M,M')$.  Moreover, $\alpha$ is the zero map (i.e., null-homotopic) if and only if it factors through $\pi(f)_*$, and hence if and only if $u$ factors through $\pi(f)$ in $\tC$.  Thus, we obtain an isomorphism $\Hom_{K(\Lambda)}(T_1,T_2) \cong \coker\ \tC(M,\pi(f)) \cong \tC(M,Y)$ by (4.4).  Furthermore, it is not hard to see that this isomorphism is compatible with the right actions of $\Endo_{K(\Lambda)}(T_1) \cong \Endo_{\tC}(M)$.  Now consider $\beta : T_2 \rightarrow T_1$ in $K(\Lambda)$.  As $\beta$ is a map of complexes, it must be induced by a unique map $v \in \tC(M',M)$ such that $v\pi(f) = 0$.  It follows that $\Hom_{K(\Lambda)}(T_2,T_1) \cong \ker\ \tC(\pi(f),M) \cong \tC(Y,M)$ by (4.3).  Again, one easily checks that this isomorphism is compatible with the left action of $\Endo_{K(\Lambda)}(T_1) \cong \Endo_{\tC}(M)$.

Presently we shall show that $\Endo_{K(\Lambda)}(T_2) \cong \Endo_{\tC}(Y)$ as rings.  First, a chain map from $T_2$ to itself is induced by a pair of maps $\alpha = \sum \alpha_i \in \tC(X,X)$ and $\beta = \sum \beta_i \in \tC(M',M')$, where $\alpha_i \in \C(X,X[i])$ and $\beta_i \in \C(M,M[i])$ are zero for almost all $i$. Furthermore, such a pair of maps commutes with the differential of $T_2$ if and only if $\beta_i f = f[i] \alpha_i$ in $\C$ for all $i \in \mathbb{Z}$.  Assuming this commutativity for a fixed $i$, and using that $\C(h,M'[i])=0$ for all $i$, we obtain a {\it unique} good map $\gamma_i$ completing a map of triangles in $\C$ by Lemma 3.1 (and the remarks preceding Lemma 3.3).
%$$\xymatrix{Y[-1] \ar[r]^h \ar[d]^{\gamma_i[-1]} & X \ar[r]^{f} \ar[d]^{\alpha_i} & M' \ar[r]^g \ar[d]^{\beta_i} & Y \ar[r] \ar[d]^{\gamma_i} & \\ Y[i-1] \ar[r]^{\pm h[i]}  & X[i] \ar[r]^{f[i]}  & M'[i] \ar[r]^{g[i]} & Y[i] \ar[r] &}$$ 
$$\xymatrix{ X \ar[r]^{f} \ar[d]^{\alpha_i} & M' \ar[r]^g \ar[d]^{\beta_i} & Y \ar[r]^h \ar[d]^{\gamma_i} & X[1]  \ar[d]^{\alpha_i[1]}  \\   X[i] \ar[r]^{f[i]}  & M'[i] \ar[r]^{g[i]} & Y[i] \ar[r]^{\pm h[i]}  & X[i+1]}$$
In particular, we see that $\alpha$ and $\beta$ induce a map $\gamma = \sum \gamma_i \in \tC(Y,Y)$, and this correspondence defines a ring homomorphism $\varphi : \Endo_{C(\Lambda)}(T_2) \rightarrow \Endo_{\tC}(Y)$.  As $g[i]$ is a right $\gen{M}$-approximation of $Y[i]$, any map $\gamma_i : Y \rightarrow Y[i]$ will lift to a map $ \beta_i : M' \rightarrow M'[i]$ making the middle square of the above diagram commute.  By Lemma 3.4, we can find a map of triangles as in the above diagram for which $\gamma_i$ is good, and this shows that $\varphi$ is surjective.  Furthermore, by Lemmas 3.2 and 3.3, the kernel of $\varphi$ consists precisely of those pairs $(\alpha,\beta)$ for which each map $(\alpha_i,\beta_i)$ is null-homotopic, i.e., those pairs that induce null-homotopic endomorphisms of the complex $T_2$.  Additionally, one can check that the isomorphisms $\Hom_{K(\Lambda)}(T_1,T_2) \cong \tC(M,Y)$  and $\Hom_{K(\Lambda)}(T_2,T_1) \cong \tC(Y,M)$ are compatible with the left and right $\Endo_{K(\Lambda)}(T_2) \cong \tC(Y,Y)$ actions, respectively, via $\varphi$. 

The proof of (2) proceeds along the same lines.  However, we need to use $M''$ in place of $M$, as it may be the case that $M' \not \in \add(M)$ in $\C$ and hence $\C(M \oplus X, M')$ may fail to be projective as an $\Endo_{\C}(M \oplus X)$-module.  (In $\tC$, on the other hand, we do have $M' \in \add(M)$.)  In place of (4.3) and (4.4), we will use the short exact sequences 
\begin{eqnarray} \ses{\C(Y,M'')}{\C(M',M'')}{\C(X,M'')}{\C(g,M'')}{\C(f,M'')} \end{eqnarray} and 
\begin{eqnarray} \ses{\C(M'',X)}{\C(M'',M')}{\C(M'',Y)}{\C(M'',f)}{\C(M'',g)}, \end{eqnarray}
which follow from (4.1) and (4.2) respectively, since $M'' \in \gen{M}$.  As before $T'$ generates $K^b(\proj \Lambda')$, as $\C(M''\oplus X,X)[1]$ is homotopic to the mapping cone of the obvious map from $\C(M''\oplus X, M') \in \add(T_1')$ to $T_2'$.  Next, notice that any map $\alpha : T_2' \rightarrow T_1'[1]$ will be induced by a map in $\C(X,M'')$ and such a map factors through $f$ by (a), making $\alpha$ null-homotopic.  Similarly, we see that $\Hom_{K(\Lambda')}(T_2',T_2'[1])=0$.  Likewise, any map $\beta : T'[1] \rightarrow T_2'$ will be induced by a map $\beta$ in $\C(M'', X)$ or $\C(M',X)$ such that $f \beta= 0$, and it follows from (4.6) that $\beta = 0$.  

To verify (ii), notice that $\Endo_{\Lambda'}(T_1') \cong \Endo_{\C}(M'')$ as rings.  Next, observe that any map $\alpha$ in $K(\Lambda')$ from $T_1'$ to $T_2'$ is induced by a unique $u \in \C(M'',M')$.  Moreover, $\alpha$ is the zero map (i.e., null-homotopic) if and only if it factors through $f_*$, and hence if and only if $u$ factors through $f$ in $\C$.  Thus, we obtain an isomorphism $\Hom_{K(\Lambda')}(T_1',T_2') \cong \coker\ \C(M'',f) \cong \C(M'',Y)$ by (4.6).  Furthermore, it is not hard to see that this isomorphism is compatible with the right actions of $\Endo_{K(\Lambda')}(T_1') \cong \Endo_{\C}(M'')$.  Now consider $\beta : T_2' \rightarrow T_1'$ in $K(\Lambda')$.  As $\beta$ is a map of complexes, it must be induced by a unique map $v \in \C(M',M'')$ such that $vf = 0$.  It follows that $\Hom_{K(\Lambda')}(T_2',T_1') \cong \ker\ \C(f,M'') \cong \C(Y,M'')$ by (4.5).  Again, one easily checks that this isomorphism is compatible with the left action of $\Endo_{K(\Lambda')}(T_1') \cong \Endo_{\C}(M'')$.  To prove that $\Endo_{K(\Lambda')}(T_2') \cong \Endo_{\C}(Y)$ as rings,   observe that  a chain map from $T_2'$ to itself is induced by a pair of maps $\alpha \in \C(X,X)$ and $\beta \in \C(M',M')$, where $\beta f = f \alpha$.   Since $\C(h,M')=0$ by (4.1), we obtain a {\it unique} good map $\gamma \in \C(Y,Y)$ completing a map of triangles in $\C$ by Lemma 3.1 (and the remarks preceding Lemma 3.3).
%$$\xymatrix{Y[-1] \ar[r]^h \ar[d]^{\gamma_i[-1]} & X \ar[r]^{f} \ar[d]^{\alpha_i} & M' \ar[r]^g \ar[d]^{\beta_i} & Y \ar[r] \ar[d]^{\gamma_i} & \\ Y[i-1] \ar[r]^{\pm h[i]}  & X[i] \ar[r]^{f[i]}  & M'[i] \ar[r]^{g[i]} & Y[i] \ar[r] &}$$ 
$$\xymatrix{ X \ar[r]^{f} \ar[d]^{\alpha} & M' \ar[r]^g \ar[d]^{\beta} & Y \ar[r]^h \ar[d]^{\gamma} & X[1]  \ar[d]^{\alpha[1]}  \\   X \ar[r]^{f}  & M' \ar[r]^{g} & Y\ar[r]^{ h}  & X[1]}$$
As above, this correspondence defines a ring homomorphism $\varphi : \Endo_{C(\Lambda')}(T_2') \rightarrow \Endo_{\C}(Y)$.  Furthermore, since $g$ is a right $\gen{M}$-approximation of $Y$, any map $\gamma: Y \rightarrow Y$ will lift to a map $ \beta: M' \rightarrow M'$ making the middle square of the above diagram commute.  By Lemma 3.4, we can find a map of triangles as in the above diagram for which $\gamma$ is good, and this shows that $\varphi$ is surjective.  By Lemmas 3.2 and 3.3, the kernel of $\varphi$ consists precisely of those pairs $(\alpha,\beta)$ that induce null-homotopic endomorphisms of the complex $T_2'$.  Additionally, one can check that the isomorphisms $\Hom_{K(\Lambda')}(T_1',T_2') \cong \C(M'',Y)$  and $\Hom_{K(\Lambda')}(T_2',T_1') \cong \C(Y,M'')$ are compatible with the left and right $\Endo_{K(\Lambda')}(T_2') \cong \C(Y,Y)$ actions, respectively, via $\varphi$.  $\Box$ \\

\noindent
{\bf Remarks.} (1) The need to consider $\gen{M}$-approximations instead of ordinary $\add(M)$-approximations can be seen in the exactness of (4.1-6) above.  Without this assumption there is no guarantee that the long exact hom sequences obtained by applying $\C(-,M)$ to the given triangle can be split up into short exact sequences.  Furthermore, it seems easier to produce examples of triangles where both maps are $\gen{M}$-approximations:  for, in case the endomorphism ring $\Lambda$ is weakly symmetric, we need only check that $f$ is a $\gen{M}$-approximation.

\vspace{1mm}
(2) Notice that the tilting complex constructed in part (1) will be isomorphic to a tilting module (in $D^b(\Lambda)$) if and only if $\tC(X,\pi(f))$ is injective, while part (2) yields a tilting module if and only if $\C(X,f)$ is injective.  In particular, $\C(X,f)$ (resp. $\tC(X,\pi(f))$) is injective if the given triangle is an AR-triangle and $X$ is not a summand of $M$ (resp. of any $M[i]$).

\vspace{1mm}
(3) If $X \stackrel{f}{\longrightarrow} M' \stackrel{g}{\longrightarrow} Y \stackrel{h}{\longrightarrow} X[1]$ is an AR-triangle, even with the added condition $X[1] \not \in \add(M' \oplus Y)$ as in Hu's and Xi's Proposition 5.1 in \cite{HuXi1}, it does not necessarily follow that $f$ and $g$ are $\gen{M'}$-approximations, as it may be the case that $X[i] \in \add(M')$ for some $i \in \mathbb{Z}$.  Hence this proposition is not a consequence of our theorem.  Nevertheless, our hypotheses (a) and (b) on the triangle $X \stackrel{f}{\longrightarrow} M' \stackrel{g}{\longrightarrow} Y \stackrel{h}{\longrightarrow} X[1]$ are typically much weaker than the assumption that it is an AR-triangle.  This is illustrated, for instance, by our examples in the following sections, where such approximation triangles can be constructed starting from any objects $X$ and $M$ of $\C$.

%(2) The only obstruction to proving the above theorem for an arbitrary triangulated category $\C$ arises in the proof that $\Endo_{C(\Lambda)}(T_2) \cong \Endo_{\tC}(Y)$.  In the general case, it can be shown that these rings contain isomorphic ideals with respect to which the factor rings are isomorphic.  It thus seems very plausible that the assumption that $\C$ is algebraic may not be truly necessary.

\vspace{1mm}
(4) The endomorphism rings of the form $\Endo_{\tC}(X \oplus M)$, taken in the orbit category, are closely related to the $\Phi$-Auslander-Yoneda algebras considered by Hu and Xi in \cite{HuXi2}, and they are a special case of the {\it perforated Yoneda algebras} studied by Hu, Koenig and Xi \cite{HuKX}.  The latter algebras, constructed with respect to an auto-equivalence $F$ of a triangulated category $\C$, are certain graded subquotients of endomorphism rings in the orbit category $\C/F$.   Moreover, given a triangle as in the above theorem, where $f$ and $g$ are  {\it cohomological approximations}, Hu, Koenig and Xi give necessary conditions for the existence of a derived equivalence between (certain factor rings of) the perforated Yoneda algebras associated to $X \oplus M$ and $Y \oplus M$ \cite{HuKX}.  

%\begin{coro}  Let $T^{\bullet} = T_1 \oplus T_2$ be a tilting complex in $K^b(\proj \Lambda)$.  Suppose that $f : T_1 \rightarrow T'$ is a left $\add(T_2)$-approximation such that the induced map $g : T' \rightarrow C(f)$ is a right $\add(T_2)$-approximation of the mapping cone $C(f)$.   Then $T_2 \oplus C(f)$ is another tilting complex in $K^b(\proj \Lambda)$. \end{coro}

\section{Applications to symmetric algebras}
\setcounter{equation}{0}

In this section, we discuss several settings where the hypotheses of Theorem 4.1 are satisfied rather easily, and we obtain interesting pairs of derived equivalent symmetric algebras.  We first focus on the case where $\C = K^b(\proj A)$ for a symmetric $k$-algebra $A$.  We start by reviewing some properties of the homotopy category of perfect complexes.

For any finite-dimensional $k$-algebra $A$, let $\nu = -\otimes_A DA$ be the Nakayama functor.  It induces an equivalence $\proj A \rightarrow \inj A$, and thus extends to an equivalence of triangulated categories $\nu : K^b(\proj A) \rightarrow K^b(\inj A)$.  Of course, if $A$ is self-injective, then $\inj A = \proj A$ and $\nu$ gives an auto-equivalence of $K^b(\proj A)$.  Furthermore, the natural isomorphism of functors $D\Hom_A(P,-) \cong \Hom_A(-,\nu P)$ for $P$ projective, also extends to a natural isomorphism $D\Hom_{K(A)}(P^{\bullet},-) \cong \Hom_{K(A)}(-,\nu P^{\bullet})$ for all $P^{\bullet} \in K^b(\proj A)$ \cite{TCRTA}.  In particular, when $A$ is self-injective $\nu$ is a Serre functor on $K^b(\proj A)$, and when $A$ is symmetric $\nu$ is the identity and we have natural isomorphisms $$D\Hom_{K(A)}(X^{\bullet},Y^{\bullet}) \cong \Hom_{K(A)}(Y^{\bullet}, X^{\bullet})$$ for all $X^{\bullet}, Y^{\bullet} \in K^b(\proj A)$  (i.e., $K^b(\proj A)$ is $0$-Calabi-Yau).

It is well-known that $K^b(\proj A)$ is an algebraic triangulated category in which the Krull-Schmidt theorem holds.  We frequently abbreviate $K^b(\proj A)$ and $\widetilde{K^b(\proj A)}$ as $K(A)$ and $\tK$ respectively.  We also point out that for any $X^{\bullet}, Y^{\bullet} \in K^b(\proj A)$, the Hom-space $\Hom_{\tK}(X,Y) = \oplus_{i \in \mathbb{Z}} \Hom_{K(A)}(X,Y[i])$ is finite-dimensional.  In particular, a left $\add(Y)$-approximation of $X$ in $\tK$ can always be constructed by summing all the maps in a basis for $\Hom_{\tK}(X,Y)$.

\begin{propos} Let $A$ be a finite-dimensional symmetric $k$-algebra.  Then for any $X^{\bullet} \in K^b(\proj A)$, the endomorphism rings $\Endo_{K(A)}(X^{\bullet})$ and $\Endo_{\tK}(X^{\bullet})$ are finite-dimensional symmetric $k$-algebras.
\end{propos}

\noindent
{\it Proof.} For the first, the $0$-Calabi-Yau property gives us $\Hom_{K(A)}(X,X) \cong D\Hom_{K(A)}(X,X)$, and the naturality in both variables implies that this is an isomorphism of $\Endo_{K(A)}(X)$-bimodules.  For the second, we have $\Hom_{K(A)}(X,X[i]) \cong D\Hom_{K(A)}(X[i],X) \cong D\Hom_{K(A)}(X,X[-i])$ for each $i \in \mathbb{Z}$.  Adding these up over all $i \in \mathbb{Z}$ yields $\Hom_{\tK}(X) \cong D\Hom_{\tK}(X)$, and again the naturality of these isomorphisms together with the definition of multiplication in $\Endo_{\tK}(X)$ shows that this is an isomorphism of bimodules. $\Box$\\

\noindent
{\bf Remarks.} (1) If $A$ is self-injective, then the Nakayama functor $\nu$ is not necessarily isomorphic to the identity.  However, we still have isomorphisms $D\Hom(X^{\bullet},Y^{\bullet}) \cong \Hom(Y^{\bullet}, \nu X^{\bullet})$, and it follows that the endomorphism rings $\Endo_{\tK}(X)$ and $\Endo_{K(A)}(X)$ are again self-injective for any $X$ with $\nu X \cong X$.  Furthermore, such an isomorphism will induce the Nakayama automorphism on the endomorphism ring.  

(2) More generally, if $\C$ is a Hom-finite triangulated category with Serre functor $\nu$, then $\Endo_{\C}(X)$ is self-injective for any $X$ with $X \cong \nu X$.  If $\nu \cong Id_{\C}$, i.e., if $\C$ is $0$-Calabi-Yau, then $\Endo_{\C}(X)$ is symmetric for any $X$. One has to be a little more careful with the total endomorphism rings in general, since these need not be finite-dimensional.\\

We now restate our main theorem from the previous section for the case where $\mathcal{C} = K^b(\proj A)$ for $A$ symmetric, as this is the case of greatest interest to us.  The proof is immediate from Theorem 4.1 and the discussion above.

\begin{therm}  Let $A$ be symmetric and let $X^{\bullet}$ and  $M^{\bullet}$ be any complexes in $K^b(\proj A)$.  Then there exists a left $\gen{M}$-approximation $f : X \rightarrow M'$ of $X$ in $K(A)$.  If $Y = C(f)$ is the mapping cone of $f$, then
 \begin{enumerate}
\item $\Endo_{\tK}(X \oplus M)$ and $\Endo_{\tK}(Y \oplus M)$ are derived equivalent symmetric algebras.
\item $\Endo_{K(A)}(X \oplus M'')$ and $\Endo_{K(A)}(Y \oplus M'')$ are derived equivalent symmetric algebras, for any $M'' \in \gen{M}$ with $M'\in\add(M'')$.
\end{enumerate}
\end{therm}

Assuming that every summand of $M$ appears in $M''$ up to a shift, the algebras in (1) are Morita equivalent to $\Endo_{\tK}(X \oplus M'')$ and $\Endo_{\tK}(Y \oplus M'')$ respectively.  The latter algebras are now derived equivalent $\mathbb{Z}$-graded algebras whose degree-$0$ subalgebras are the algebras in (2), which are again derived equivalent.  We find it interesting that the degree-$0$ subalgebras will typically have more simples (up to isomorphism) than the full graded endomorphism rings.  This is because $M'$ (and hence $M''$) may contain multiple shifts of the same complex as summands, yielding summands that are isomorphic in $\tK$ but not in $K(A)$.  

\vspace{3mm}
In looking for applications of our theorem we do not need to restrict our attention to categories of perfect complexes.  In fact, Happel has shown that $K^-(\proj A)$ is equivalent to a full triangulated subcategory of the stable category of the repetitive algebra of $A$, which is also equivalent to the stable category of $\mathbb{Z}$-graded modules over the trivial extension $T(A) = A \oplus DA$, where $A$ is in degree $0$ and $DA$ is in degree $1$ \cite{TCRTA}.  Thus the category of perfect complexes can be viewed as a triangulated subcategory of this stable category.  %We feel however that the nature of the suspension functors makes computations of approximations and endomorphism rings simpler in categories of complexes than in stable categories.
This suggests that it may be interesting to apply Theorem 4.1 to stable categories of $\mathbb{Z}$-graded modules over (nice) self-injective algebras more generally.  In such contexts the $\mathbb{Z}$-grading will often guarantee existence of the necessary approximations as well as the finite-dimensionality of the endomorphism rings in the orbit category: for instance, it suffices to know that for any graded modules $X, Y$ we have $X$ and $\Omega^i Y$ concentrated in disjoint degrees for all $|i|$ sufficiently large.

On the other hand, if we wish to apply Theorem 4.1 to the (ungraded) stable category $\stmod A$ of a self-injective algebra $A$, the necessary approximations in the orbit category may fail to exist.  Even if they do exist, the endomorphism rings in the orbit category will have the form $\oplus_{i \in \mathbb{Z}}\stHom_A(\Omega^iX,X)$, and will most likely not be finite-dimensional.  However, if we assume that $M$ is $\Omega$-periodic, then there is no problem obtaining a left $\gen{M}$-approximation of any $X$ in $\stmod A$.  Then we can still apply Theorem 4.1(2) starting with $M$ and any $A$-module $X$.  The endomorphism ring $\Lambda' = \stEnd_A(M' \oplus X)$, for $M' \in \gen{M}$, will be weakly symmetric provided the Serre functor $\tau \Omega^{-1} \cong \nu \Omega$ fixes each indecomposable summand of $M \oplus X$ (since the Serre functor is a triangulated functor, it follows that it also fixes any object in $\gen{M}$).  We now summarize these observations in a corollary.

\begin{coro} Let $X$ and $M$ be modules over a self-injective algebra $A$ with $\Omega^n M \cong M$ for some $n \geq 1$ and $\nu \Omega Z \cong Z$ for each indecomposable summand of $X \oplus M$.  Let $\tilde{M} = \oplus_{i=0}^{n-1} \Omega^iM$ and suppose $f : X \rightarrow M'$ is a left $\add(\tilde{M})$-approximation of $X$ in $\stmod A$ and $Y$ is the cone of $f$, then $\stEnd_A(X \oplus M'')$ and $\stEnd_A(Y \oplus M'')$ are derived equivalent weakly symmetric algebras for any $M'' \in \add(\tilde{M})$ with $M' \in \add(M'')$.
\end{coro}

%\noindent {\bf Remark.}  We point out that $Y$ can also be obtained as the cokernel of a left $\add(A \oplus M)$-approximation $\tilde{f}$ of $X$ in $\rmod A$.

While the hypotheses in the above corollary appear to be somewhat restrictive in $\stmod A$, they turn out to be satisfied automatically in some categories of Cohen-Macaulay modules.  Namely, let $R$ be an isolated Gorenstein hypersurface singularity of dimension $d$, and let $\C = \und{CM}(R)$ be the stable category of (maximal) Cohen-Macaulay $R$-modules.  When $d$ is odd, $\C$ is a Hom-finite $0$-Calabi-Yau triangulated category with suspension $\Sigma = \Omega^{-1}$ satisfying $\Sigma^2 \cong Id$ (see for instance \cite{BIKR}).  The same argument as above now gives the following.

\begin{coro}  Let $\C = \und{CM}(R)$ for an odd-dimensional isolated Gorenstein hypersurface singularity $R$.  For any $X, M \in \C$, there exists a left $\add(M \oplus \Omega M)$-approximation $f : X \rightarrow M'$ of $X$ in $\C$.  If $Y$ is the cone of $f$ then $\stEnd_R(X \oplus M'')$ and $\stEnd_R(Y \oplus M'')$ are derived equivalent symmetric algebras for any $M'' \in \add(M \oplus \Omega M)$ with $M' \in \add(M'')$.
\end{coro}

In \cite{BIKR}, Burban, Iyama, Keller and Reiten consider endomorphism rings of cluster-tilting objects in $\und{CM}(R)$ when $R$ is a curve singularity.  Recall that $T \in CM(R)$ is a cluster-tilting object if 
$$\add(T)  =  \{ X \in CM(R)\ |\ \Ext^1_R(X,T)=0\} = \{Y \in CM(R)\ |\ \Ext^1_R(T,Y)=0\}.$$ 
When $R$ is an isolated Gorenstein singularity of dimension $d \leq 3$, Iyama \cite{Iyama2} has shown that the endomorphism rings of cluster-tilting objects in $CM(R)$ are all derived equivalent (as in Hu's and Xi's Theorem, these derived equivalences are furnished by tilting modules of projective dimension $1$).  In this context, our corollary establishes the derived equivalence of the {\it stable} endomorphism rings of cluster-tilting objects which are connected by a mutation (as in Definition 1.2 of \cite{BIKR}).  

\begin{coro}  Let $R$ be an odd-dimensional Gorenstein hypersurface that is an isolated singularity, and let $T=M \oplus X$ be a basic cluster-tilting object in $CM(R)$ with $R \in \add(M)$ and $X$ indecomposable.  Let $Y$ denote the cokernel of a minimal left $\add(M)$-approximation $f : X \rightarrow M'$ of $X$ in $CM(R)$.  Then $T' = M \oplus Y$ is a cluster-tilting object in $CM(R)$ and $\stEnd_R(T')$ is derived equivalent to $\stEnd_R(T)$.
\end{coro}

\noindent
{\it Proof.}  Let $\C = \und{CM}(R)$.  Since $T$ is a cluster-tilting object and $T[2] \cong T$ in $\C$, we have $\C(T,T[-1]) = 0$.  Thus $\ul{f}$ is a left $\add(M \oplus \Omega M)$-approximation of $X$ in $\C$, and its cone is isomorphic to $Y$.  Now apply the previous corollary with $M''=M$.  $\Box$ \\

  If $T$ and $T'$ are two cluster-tilting objects in $CM(R)$, the tilting module yielding the derived equivalence between $\Endo_R(T)$ and $\Endo_R(T')$ is $\Hom_R(T,T')$.  When $T$ and $T'$ are related by a mutation, the tilting complex yielding the derived equivalence between $\stEnd_R(T)$ and $\stEnd_R(T')$ (as in Theorem 4.1(2)) is given by $P^{\bullet} \otimes_{\Endo_R(T)} \stEnd_R(T)$, where $P^{\bullet}$ is the projective resolution of $\Hom_R(T,T')$ over $\Endo_R(T)$.  When $T$ and $T'$ are not connected by a single mutation, we do not know whether the same construction still yields a tilting complex over $\stEnd_R(T)$ with endomorphism ring isomorphic to $\stEnd_R(T')$.

\section{Examples}
\setcounter{equation}{0}

One nice feature of Theorem 5.2 is that it can be applied to {\it any} complexes of projectives over {\it any} symmetric algebra.  In fact we will see here that even starting from relatively simple local, commutative algebras $A$, we are able to produce many interesting examples of derived equivalent symmetric algebras, generalizing some known families.  In particular, the derived equivalences in Examples 1 and 2 with $m\neq 1, n-1$ appear to be new.

\vspace{3mm}
\noindent
{\bf Example 1.} Let $A = k[x]/(x^n)$ for some $n \geq 2$, and set $T = T(m) = \ul{A} \oplus [A \stackrel{x^m}{\longrightarrow} \ul{A}] \in K^b(\proj A)$ for any $m$ with $1 \leq m < n$.  We consider the following maps between the summands of $T$ in $\add_{\tK}(T)$.

$$\begin{array}{ccccc} \xymatrix{A \save+<-4ex,-4ex> \drop{\alpha:} \restore \ar[d]^x \\ A}\ & \xymatrix{0 \save+<-4ex,-4ex> \drop{\beta:} \restore \ar[r] \ar[d] & A \ar[d]^1 \\ A \ar[r]^{x^m} & A}\ & \xymatrix{A \save+<-4ex,-4ex> \drop{\gamma:} \restore \ar[d]^1 \ar[r]^{x^m} & A \ar[d] \\ A \ar[r] & 0}\ &  \xymatrix{A \save+<-4ex,-4ex> \drop{\delta:} \restore \ar[d]^x \ar[r]^{x^m} & A \ar[d]^x \\ A \ar[r]_{x^m} & A}\ &  \xymatrix{A \save+<-4ex,-4ex> \drop{\epsilon:} \restore \ar[d] \ar[r]^{x^m} & A \ar[r] \ar[d]^{x^{n-m}} & 0 \ar[d] \\ 0 \ar[r] & A \ar[r]_{x^m} & A}
\end{array}$$

It is not hard to see that these generate $\Endo_{\tK}(T)$ and are irreducible (except in some degenerate cases--see below).  We thus see that $\Lambda(m) := \Endo_{\tK}(T)$ is given by the quiver and relations:

$$\xymatrixcolsep{1.0pc} \xymatrix{ 1 \ar@(dl,ul)^{\alpha} \ar[rr]<0.5ex>^{\beta} & & 2  \ar[ll]<0.5ex>^{\gamma}  \ar@(r,u)_{\delta} \ar@(r,d)^{\epsilon} }$$

\begin{eqnarray*} \beta \alpha^m = \gamma \beta = \delta^m = \epsilon^2 = 0,\\
\beta \alpha = \delta \beta,\ \delta \epsilon = \epsilon \delta,\ \gamma \delta = \alpha \gamma,\\\alpha^{n-m} = \gamma \epsilon \beta,\ \delta^{n-m} = \beta \gamma \epsilon + \epsilon \beta \gamma
\end{eqnarray*}

Writing $T_1 = \ul{A}$ and $T_2 = [A \stackrel{x^m}{\longrightarrow} \ul{A}] $ for the summands of $T$, we find that $u := \bnc{\gamma}{\gamma\epsilon}  : T_2 \rightarrow T_1[1] \oplus T_1$ is a left $\gen{T_1}$-approximation, yielding the following triangle in $K^b(\proj A)$.

\begin{eqnarray}% \vcenter{\xymatrix{ A \ar[d]^{x^m} \ar[r]^1 & A \ar[d]^{0} \ar[r]^{x^m} & A \ar[d]^{x^{n-m}} \ar[r] & \\ A \ar[r]_{x^{n-m}} & A \ar[r]_1 & A \ar[r] &}} 
(A \stackrel{x^m}{\rightarrow} \ul{A} ) \stackrel{(1,x^{n-m})}{\longrightarrow} A \oplus A[1] \stackrel{(x^m,1)}{\longrightarrow} (A \stackrel{x^{n-m}}{\rightarrow} \ul{A}) \longrightarrow 
\end{eqnarray}

Theorem 5.2(1) thus shows that $\Lambda(m)$ and $\Lambda(n-m)$ are derived equivalent.  On the other hand, the left $\add_{\tK}(T_2)$-approximation of $T_1$ is given by $\beta$, which has mapping cone isomorphic to $T_1[1]$.  Thus the two total endomorphism rings in this case will be isomorphic, and we obtain a nontrivial auto-equivalence of the corresponding derived category.

The degenerate cases alluded to above occur for $m =1, n-1$.  When $m=1$, we have $\delta^1=0$, and hence $\Lambda(1)$ coincides with the {\it dihedral algebra} $D(2\B)^{n-1,0}$ as introduced in \cite{Erd}.  When $m=n-1$, we have $\alpha^1 = \gamma \epsilon \beta$ and $\delta^1 = \beta \gamma \epsilon + \epsilon \beta \gamma$.  In addition, $\delta^{n-1} = 0 = \epsilon^2 = \gamma \beta$ now imply that $(\beta \gamma \epsilon)^{n-1} + (\epsilon \beta \gamma)^{n-1} = 0$, and hence $\Lambda(n-1)$ is isomorphic to the dihedral algebra $D(2\A)^{n-1,0}$.  These algebras were first shown to be derived equivalent by Holm \cite{Holm1}.

\vspace{3mm}
\noindent
{\bf Example 2.} We now illustrate Theorem 5.2(2) in the context of the previous example.  We again use the triangle (6.1), to conclude that the ordinary endomorphism rings of $T(m) := A  \oplus [A \stackrel{x^m}{\longrightarrow} \ul{A}] \oplus A[1] $ and $T(n-m) = A \oplus [A \stackrel{x^{n-m}}{\longrightarrow} \ul{A}] \oplus A[1] $ are derived equivalent.  To describe this endomorphism ring for a fixed $m$ we define (usually) irreducible maps as follows:

$$\begin{array}{ccccccc} \xymatrix{A \save+<-4ex,-4ex> \drop{\alpha:} \restore \ar[d]^x \\ A}\ & \xymatrix{A[1] \save+<-4ex,-4ex> \drop{\alpha':} \restore \ar[d]^x \\ A[1]}\ & \xymatrixcolsep{1.5pc} \xymatrix{0 \save+<-4ex,-4ex> \drop{\beta:} \restore \ar[r] \ar[d] & A \ar[d]^1 \\ A \ar[r]^{x^m} & A}\ & \xymatrixcolsep{1.5pc} \xymatrix{A \save+<-4ex,-4ex> \drop{\beta':} \restore \ar[r] \ar[d]^{x^{n-m}} & 0 \ar[d] \\ A \ar[r]_{x^m} & A}\ & \xymatrixcolsep{1.5pc} \xymatrix{A \save+<-4ex,-4ex> \drop{\gamma:} \restore \ar[d]^1 \ar[r]^{x^m} & A \ar[d] \\ A \ar[r] & 0}\ & \xymatrixcolsep{1.5pc} \xymatrix{A \save+<-4ex,-4ex> \drop{\gamma':} \restore \ar[d] \ar[r]^{x^m} & A \ar[d]^{x^{n-m}} \\ 0 \ar[r] & A}\ & \xymatrixcolsep{1.5pc} \xymatrix{A \save+<-4ex,-4ex> \drop{\delta:} \restore \ar[d]^x \ar[r]^{x^m} & A \ar[d]^x \\ A \ar[r]_{x^m} & A}
\end{array}$$

Then $\Gamma(m) := \Endo_{K(A)}(T(m))$ is given by the quiver and relations below for each $1 \leq m \leq n-1$.  By Theorem 5.2(2), $\Gamma(m)$ is derived equivalent to $\Gamma(n-m)$ for each positive integer $m$ smaller than a fixed $n$.

$$ \xymatrix{ 1 \ar@(dl,ul)^{\alpha} \ar[r]<0.5ex>^{\beta} & 2 \ar[l]<0.5ex>^{\gamma'}  \ar@(ul,ur)^{\delta} \ar[r]<0.5ex>^{\gamma} & 3 \ar[l]<0.5ex>^{\beta'} \ar@(dr,ur)_{\alpha'} }$$

\begin{eqnarray*} \beta \alpha^m = \beta'(\alpha')^m = \alpha^m\gamma' = (\alpha')^m\gamma = \gamma \beta = \gamma'\beta' = \delta^m = 0, \\ \alpha^{n-m} = \gamma' \beta,\ (\alpha')^{n-m} = \gamma \beta',\ \delta^{n-m} = \beta \gamma' + \beta'\gamma, \\ \beta \alpha = \delta \beta, \beta'\alpha' = \delta \beta', \gamma\delta = \alpha'\gamma, \gamma'\delta = \alpha \gamma'.
\end{eqnarray*}

As before, when $m=1,n-1$ the quiver presentation of $\Gamma(m)$ is somewhat redundant.  In particular, when $m=1$, we have $\delta = 0$ and one can see that $\Gamma(1)$ is isomorphic to the algebra of dihedral type $D(3\D)_2^{1,1,n-1,n-1}$ with $3$ simples.  Similarly, when $m=n-1$, we have $\delta = \beta \gamma'+\beta'\gamma, \alpha = \gamma'\beta$ and $\alpha'=\gamma\beta'$, from which we see that $\Gamma(n-1)$ is isomorphic to the algebra of dihedral type $D(3\A)_2^{n-1,n-1}$.  These algebras are among the much larger family of (non-block) algebras of dihedral type with 3 simples that Holm shows are derived equivalent in \cite{Holm2}.

\vspace{3mm}
\noindent
{\bf Example 3.}  We now show how to realize all standard algebras of dihedral type with two simples as endomorphism rings of complexes.  Set $A = k[x,y]/(x^n-y^s, xy)$ and consider $T = \ul{A} \oplus [A \stackrel{x}{\longrightarrow} \ul{A}]$.  To describe $\Lambda(n,s) = \Endo_{\tK}(T)$ we first identify the irreducible maps between the summands of $T$ and their shifts.

$$\begin{array}{ccccc} \xymatrix{A \save+<-4ex,-4ex> \drop{\alpha:} \restore \ar[d]^x \\ A}\ & \xymatrix{0 \save+<-4ex,-4ex> \drop{\beta:} \restore \ar[r] \ar[d] & A \ar[d]^1 \\ A \ar[r]^{x} & A}\ & \xymatrix{A \save+<-4ex,-4ex> \drop{\gamma:} \restore \ar[d]^1 \ar[r]^{x} & A \ar[d] \\ A \ar[r] & 0}\ &  \xymatrix{A \save+<-4ex,-4ex> \drop{\epsilon:} \restore \ar[d] \ar[r]^{x} & A \ar[r] \ar[d]^{y} & 0 \ar[d] \\ 0 \ar[r] & A \ar[r]_{x} & A}
\end{array}$$

From these one computes the following quiver and relations for $\Lambda(n,s)$.

$$\xymatrixcolsep{1.0pc} \xymatrix{ 1 \ar@(dl,ul)^{\alpha} \ar[rr]<0.5ex>^{\beta} & & 2  \ar[ll]<0.5ex>^{\gamma} \ar@(ur,dr)^{\epsilon} }$$

$$\alpha \gamma = \gamma \beta = \beta \alpha = \epsilon^2=0, \ \ \alpha^n = (\gamma \epsilon \beta)^s, \ \   (\epsilon \beta \gamma)^s + (\beta \gamma \epsilon)^s = 0.$$

Hence $\Lambda(n,s)$ is isomorphic to the dihedral algebra $D(2\B)^{s,n}(0)$ with two simples \cite{Holm2}.  As in Example 1, writing $T_1$ and $T_2$ for the two summands of $T$, we have $u := \bnc{\gamma}{\gamma\epsilon}  : T_2 \rightarrow T_1[1] \oplus T_1$ is a left $\add_{\tK}(T_1)$-approximation, yielding the following triangle in $K^b(\proj A)$.

\begin{eqnarray}% \vcenter{ \xymatrix{ A \ar[d]^{x} \ar[r]^1 & A \ar[d]^{0} \ar[r]^{x} & A \ar[d]^{y} \ar[r] & \\ A \ar[r]_{y} & A \ar[r]_1 & A \ar[r] &}}
(A \stackrel{x}{\rightarrow} \ul{A}) \stackrel{(1,y)}{\longrightarrow} A \oplus A[1] \stackrel{(x,1)}{\longrightarrow} (A \stackrel{y}{\rightarrow} \ul{A}) \longrightarrow
 \end{eqnarray}

Theorem 5.2(1) thus shows that $\Lambda(n,s)$ is derived equivalent to $\Endo_{\tK}(A \oplus [A \stackrel{y}{\longrightarrow}\und{A}])$, which is isomorphic to $\Lambda(s,n)$ via the isomorphism $k[x,y]/(x^n-y^s, xy) \rightarrow k[x,y]/(x^s-y^n, xy)$ interchanging $x$ and $y$.  This derived equivalence appears also in \cite{Holm2}, Lemma 3.2.

\vspace{3mm}
It should be noted that it is also possible--and not any more difficult--to start with one of the symmetric algebras $\Lambda$ from the above examples, write down the appropriate Okuyama-Rickard tilting complex and compute its endomorphism ring to obtain the corresponding algebra derived equivalent to $\Lambda$.  Our approach to constructing these derived equivalences differs mainly in its point of view.  For one, it offers a useful way of deducing that two algebras--defined only abstractly as endomorphism rings (as in Corollary 5.5 for instance)--must be derived equivalent, without computing them explicitly.   Furthermore, in the examples above, we see that these derived equivalences are essentially controlled by homological properties of much simpler local commutative $k$-algebras.  By viewing the algebras $\Lambda(1)$ and $\Gamma(1)$ in Examples 1 and 2 as endomorphism rings of complexes over the ring $A$, we obtain a natural way of generalizing these algebras by modifying the complexes whose endomorphism rings we consider, and we obtain a family of derived equivalent pairs of algebras with a single computation.  As it would be interesting to understand to what extent similar results hold for derived equivalences between symmetric algebras furnished by Okuyama-Rickard complexes, we propose the following problems for further study.

\vspace{3mm}
\noindent
{\bf Problems.}  (1) Determine which finite-dimensional symmetric $k$-algebras $\Lambda$ arise as endomorphism rings of objects in $K^b(\proj A)$ or $\tilde{K}^b(\proj A)$ for a symmetric {\it local} $k$-algebra $A$.  \\
(2) If $\Lambda$ arises an endomorphism ring in this way, does every Okuyama-Rickard tilting complex for $\Lambda$ arise from some triangle in $K^b(\proj A)$ as in Theorem 4.1?

\vspace{2mm}
A necessary condition for $\Lambda \cong \Endo_{\tilde{K}(A)}(T)$ with $A$ local is that the entries of the Cartan matrix of $\Lambda$ are all at least $2$.  This is a result of the fact that for any two bounded complexes $X^{\bullet}$ and $Y^{\bullet}$ of projective $A$-modules (with radical maps for differentials) there is at least one non-null-homotopic map from the right-most term (resp., left-most term) of $X^{\bullet}$ to the left-most term (resp., right most term) of $Y^{\bullet}$.  Thus there are many symmetric algebras that cannot be realized in this way.  However, for such an algebra $\Lambda$, a positive answer to (2) could aid in the classification of the algebras that are linked to $\Lambda$ by a sequence of tilting mutations as defined in \cite{Aih1}.

\section{Interpretation via differential graded algebras}
\setcounter{equation}{0}

We conclude this article by looking at Theorem 5.2 from the point of view of differential graded algebras.  Recall that a {\bf dg-algebra} is a $\mathbb{Z}$-graded $k$-algebra $\mathbf{A} = \oplus_{i \in \mathbb{Z}} \mathbf{A}_i$ with a differential $d : \mathbf{A}_i \rightarrow \mathbf{A}_{i+1}$ satisfying $d(ab) = d(a)b + (-1)^{|a|}ad(b)$ for all homogenous elements $a,b \in \mathbf{A}$, where $|a|$ denotes the degree of $a$.  The cohomology $H^*(\mathbf{A}) = H^*(\mathbf{A},d)$ inherits the structure of a $\mathbb{Z}$-graded $k$-algebra from $\mathbf{A}$.

A key example of a dg-algebra can be constructed from any complex $(C^{\bullet}, \delta)$ of $R$-modules.  We set $\mathbf{A}_i := \prod_{n \in \mathbb{Z}} \Hom_R(C_n, C_{n+i})$ and define $d((f_n)_n) := (\delta f_n - (-1)^i f_{n+1}\delta)_n$ for all $(f_n)_n \in \mathbf{A}_i$.  It is easily checked that $(\mathbf{A},d)$ is a dg-algebra, often denoted $\RHom_R(C^{\bullet},C^{\bullet})$.  Moreover, $H^i(\mathbf{A}) \cong \Hom_{K(R)}(C, C[i])$ for all $i$.  In the notation of Section 4 we thus have $H^0(\mathbf{A}) \cong \Endo_{K(R)}(C)$ and $H^*(\mathbf{A}) = \Endo_{\tilde{K}(R)}(C)$.

In particular, Theorem 5.2 yields two different dg-algebras $\mathbf{\Lambda} = \RHom_A(X \oplus M',X \oplus M')$ and $\mathbf{\Gamma} = \RHom_A(Y \oplus M',Y \oplus M')$ whose cohomology rings are derived equivalent.  Furthermore, the degree-$0$ subalgebras of these cohomology rings are also derived equivalent.  This motivates the following problem.

\vspace{3mm}
\noindent
{\bf Problem.} Describe the tilting procedures of Theorem 5.2, yielding derived equivalences between the cohomology rings $H^*(\mathbf{\Lambda})$ and $H^*(\mathbf{\Gamma})$, on the level of the dg-algebras $\mathbf{\Lambda}$ and $\mathbf{\Gamma}$.  In particular, is the derived equivalence between the cohomology rings of $\mathbf{\Lambda}$ and $\mathbf{\Gamma}$ a shadow of some deeper type of equivalence between $\mathbf{\Lambda}$ and $\mathbf{\Gamma}$?  More generally, what relations can be imposed between two dg-algebras to ensure that their cohomology rings (or their degree-$0$ cohomology rings) are derived equivalent?\\

To illustrate these two dg-algebras in one concrete case, we turn to Example 3 from the previous section.  We have $A = k[x,y]/(x^n-y^s, xy)$.  Notice that, as $A$-modules, both complexes $$T^{(x)} = \ul{A} \oplus [A \stackrel{x}{\longrightarrow} \ul{A}]\ \mbox{and}\ T^{(y)} = \ul{A} \oplus [A \stackrel{y}{\longrightarrow} \ul{A}]$$ are isomorphic to $A^3$, and hence their endomorphism rings are naturally identified with $M_3(A)$.  Below, we indicate the degree of each component of these matrix dg-algebras using a subscript, and we indicate the action of the differential with arrows.  If there is no arrow leaving some component, then the differential vanishes on that component.

$$\begin{array}{ccc}  \xymatrixrowsep{1.0pc} \xymatrixcolsep{1.0pc} \xymatrix{ A_0 \save[0,0].[2,0]!C *\frm{(}="Vbr" \restore & A_1 & A_0 \ar[l]_{-x} \save[0,0].[2,0]!C *\frm{)}="Vbr" \restore \\ A_{-1} \ar[d]^x & A_0 \ar[d]_x & A_{-1} \ar[l]_x \ar[d]_x \\ A_0 & A_1 & A_{0} \ar[l]_{-x}} 
 & \hspace{2cm} & \xymatrixrowsep{1.0pc} \xymatrixcolsep{1.0pc} \xymatrix{ A_0 \save[0,0].[2,0]!C *\frm{(}="Vbr" \restore & A_1 & A_0 \ar[l]_{-y} \save[0,0].[2,0]!C *\frm{)}="Vbr" \restore \\ A_{-1} \ar[d]^y & A_0 \ar[d]_y & A_{-1} \ar[l]_{y} \ar[d]_y \\ A_0 & A_1 & A_{0} \ar[l]_{-y}}  \\ \\
\RHom_A(T^{(x)},T^{(x)}) & \hspace{2cm} & \RHom_A(T^{(y)},T^{(y)})
\end{array}$$
In particular, we see that these two dg-algebras are isomorphic as graded algebras, and differ only in the actions of their differentials.  Of course, we would like to know whether there is a direct construction which produces one of these from the other, and which could give an alternate explanation of the derived equivalence of their cohomology rings.

\end{document}